\documentclass [12pt] {amsart}

\usepackage [T1] {fontenc}
\usepackage {amssymb}
\usepackage {amsmath}
\usepackage {amsthm}
\usepackage{tikz}

\usepackage{url}

 \usepackage[titletoc,toc,title]{appendix}

\usetikzlibrary{arrows,decorations.pathmorphing,backgrounds,positioning,fit,petri}
\usetikzlibrary{matrix}

\usepackage{enumitem}
\usepackage{hyperref}

\DeclareMathOperator {\td} {td}

\DeclareMathOperator {\im} {Im}
\DeclareMathOperator {\pr} {pr}
\DeclareMathOperator {\ppr} {Pr}

\DeclareMathOperator {\Der} {Der}

\DeclareMathOperator {\SL} {SL}

\DeclareMathOperator {\GL} {GL}

\DeclareMathOperator {\Ann} {Ann}

\DeclareMathOperator {\Exp} {Exp}

\DeclareMathOperator {\alg} {alg}

\DeclareMathOperator {\Jac} {Jac}
\DeclareMathOperator {\rk} {rk}

\DeclareMathOperator {\Span} {span}
\DeclareMathOperator {\C} {\mathbb{C}}

\DeclareMathOperator {\h} {\mathbb{H}}
\DeclareMathOperator {\Z} {\mathbb{Z}}
\DeclareMathOperator {\Q} {\mathbb{Q}}

\DeclareMathOperator {\G} {\mathbb{G}}
\DeclareMathOperator {\aaa} {a}
\DeclareMathOperator {\m} {m}

\DeclareMathOperator {\rar} {\rightarrow}
\DeclareMathOperator {\seq} {\subseteq}

\DeclareMathAlphabet\urwscr{U}{urwchancal}{m}{n}%
\DeclareMathAlphabet\rsfscr{U}{rsfso}{m}{n}
\DeclareMathAlphabet\euscr{U}{eus}{m}{n}
\DeclareFontEncoding{LS2}{}{}
\DeclareFontSubstitution{LS2}{stix}{m}{n}
\DeclareMathAlphabet\stixcal{LS2}{stixcal}{m} {n}

\theoremstyle {definition}
\newtheorem {definition}{Definition} [section]

\newtheorem* {claim} {Claim}

\theoremstyle {plain}

\newtheorem {lemma} [definition] {Lemma}
\newtheorem {theorem} [definition] {Theorem}
\newtheorem {fact} [definition] {Fact}

\theoremstyle {remark}
\newtheorem {remark} [definition] {Remark}

\calclayout

\begin {document}

\title{Differential Existential Closedness for the $j$-function}

\author{Vahagn Aslanyan}
\address{Vahagn Aslanyan, School of Mathematics, University of East Anglia, Norwich, NR4 7TJ, UK}
\email{V.Aslanyan@uea.ac.uk}
\author{Sebastian Eterovi\'c}
\address{Sebastian Eterovi\'c, Department of Mathematics, UC Berkeley, CA 94720, USA}
\email{eterovic@berkeley.edu}
\author{Jonathan Kirby}
\address{Jonathan Kirby, School of Mathematics, University of East Anglia, Norwich, NR4 7TJ, UK}
\email{Jonathan.Kirby@uea.ac.uk}

\thanks{VA and JK were supported by EPSRC grant EP/S017313/1. SE was partially supported by NSF RTG grant DMS-1646385}

\date{\today}



%\vspace*{-1cm}

\maketitle

\begin{abstract}
We prove the Existential Closedness conjecture for the differential equation of the $j$-function and its derivatives. It states that in a differentially closed field certain equations involving the differential equation of the $j$-function have solutions. Its consequences include a complete axiomatisation of $j$-reducts of differentially closed fields, a dichotomy result for strongly minimal sets in those reducts, and a functional analogue of the Modular Zilber-Pink with Derivatives conjecture.
\end{abstract}

\section{Introduction}

Let $\h$ be the complex upper half-plane and $j:\h \rightarrow \C$ be the modular $j$-function. The transcendence properties of this function have been well studied. For example, it is well known that for a matrix $g\in \GL_2(\C)$ the functions $j(z)$ and $j(gz)$ are algebraically dependent if and only if $g$ has rational entries (and positive determinant), and in that case the dependence is given by a \emph{modular polynomial}. It is also known that $j(z), j'(z), j''(z)$ are algebraically independent functions over $\C(z)$ and $j'''\in \Q(j,j',j'')$ \cite{mahler}. These statements are generalised by a theorem of Pila and Tsimerman \cite{Pila-Tsim-Ax-j}, known as Ax-Schanuel for the $j$-function, for it is an analogue of the Ax-Schanuel theorem for the exponential function \cite{Ax}. A precise statement of Ax-Schanuel will be given in Section \ref{section-ax-sch}, but we give a rough description now. The $j$-function satisfies an order three differential equation (see Section \ref{section-ax-sch}), and the Ax-Schanuel theorem is a transcendence statement about its solutions often formulated in a differential algebraic language. It can be seen as a necessary condition for a system of differential equations in terms of the equation of the $j$-function to have a solution in a differential field $F$. In particular, it implies that ``overdetermined'' systems cannot have a solution.

Let $E_j(x,y)$ denote the differential equation of the $j$-function. Here we think of $x$ as some complex function and of $y$ as $j(x)$, although a solution to the differential equation could be any function $j(gx)$ for any $g\in \GL_2(\C)$. Solving a system of equations is equivalent to finding points of the form $\left(\bar{z},\bar{j}\right)$ in an algebraic variety $V\subseteq F^{2n}$ such that each $(z_i,j_i)$ is in $E_j$.  The notion of an overdetermined system can be captured by some dimension conditions on the algebraic variety $V\seq F^{2n}$ dictated by the Ax-Schanuel theorem, and varieties not satisfying those conditions are called \emph{$j$-broad}. $j$-broadness of $V$ means that $V$ is not too small, in particular, it implies that $\dim V\geq n$. See Section \ref{sect-EC-j-noder} for a definition.

Thus, in some sense Ax-Schanuel states that if $V$ contains a ``sufficiently generic'' $E_j$-point, then $V$ must be $j$-broad. There is a dual statement according to which $j$-broad varieties do contain $E_j$-points in suitably large differential fields, namely, differentially closed fields. This dual statement is known as \emph{Existential Closedness}, henceforth referred to as EC, and was conjectured in \cite{Aslanyan-adequate-predim}.\footnote{Structures satisfying EC may not be existentially closed (in a first-order sense) in arbitrary extensions, but they are existentially closed in so called ``strong'' extensions. See \cite{Aslanyan-adequate-predim,Kirby-semiab} for details.} In this paper we prove it and establish several related results. The following is one of the main theorems of the paper. 

\begin{theorem}\label{thm-intro-ec-j}
Let $F$ be a differential field,\footnote{In this paper all fields are of characteristic $0$. } and $V\subseteq F^{2n}$ be a $j$-broad variety. Then there is a differential field extension $K\supseteq F$ such that $V(K)$ contains an $E_j$-point. In particular, if $F$ is differentially closed then $V(F)$ contains an $E_j$-point.
\end{theorem}

%\footnote{ {When we write $V \seq F^{2n}$, we mean that $V$ is defined over $F$ and is identified with its $F^{\alg}$-points.}} interpreted in a differential field as the set of all $4$-tuples $(z,j,j',j'')$ where $(z,j)$ is an $E_j$-point and $j', j''$ are thought of as the first and second derivatives of $j$ with respect to $z$.\footnote{Here $'$ is not used to denote a derivative of a function, and $y',y''$ denote coordinates/variables, while $j',j''$ are some elements of a field. Thus, $j'$ is not the derivative of $j$, but it is thought of as the derivative of $j$ by $z$ in an abstract differential field, which will be defined rigorously in Section \ref{section-ax-sch}.}

We actually prove a more general theorem incorporating the $j$-function and its first two derivatives. Let $E_J(x,y,y',y'')$ be a $4$-ary relation where, by definition, $(x,y)\in E_j$ and we think of $y'$ and $y''$ as the first and second derivatives of $y$ with respect to $x$. In other words, we may think of $x$ as a complex function, and of $y,y',y''$ as $j(x),j'(x),j''(x)$ where $j', j''$ are the first two derivatives of the $j$-function. We do not consider further derivatives since $j'''$ is already algebraic over $j,j',j''$. The Ax-Schanuel theorem is in fact a transcendence statement for $E_J$-points in a differential field. As above, we define $J$-\emph{broad} varieties $V\subseteq F^{4n}$ by some dimension conditions dictated by Ax-Schanuel.\footnote{The use of upper case $J$ in some notations and terms is motivated by the fact that the vector function $(j,j',j''):\h \rar \C^3$ is often denoted by $J$.} In particular, if $V\seq F^{4n}$ is $J$-broad then $\dim V \geq 3n+1$. See Section \ref{sect-broad-free-def} for a precise definition.\footnote{Note that both $j$-broadness and $J$-broadness are algebraic (as opposed to differential) properties of a variety $V$ and they just state that certain projections of $V$ have sufficiently large dimension.} Then we have the following EC statement for $E_J$ which was conjectured in \cite{Aslanyan-adequate-predim}.

\begin{theorem}\label{thm-intro-ec-J}
Let $F$ be a differential field, and $V\subseteq F^{4n}$ be a $J$-broad variety. Then there is a differential field extension $K\supseteq F$ such that $V(K)$ contains an $E_J$-point. In particular, if $F$ is differentially closed then $V(F)$ contains an $E_J$-point.
\end{theorem}

Further, we prove that certain varieties contain non-constant $E_J$-points. Those varieties are called \emph{strongly $J$-broad}. Such a variety (in $F^{4n}$) must have dimension $\geq 3n+1$.

\begin{theorem}\label{thm-intro-ec-J-nonconst}
Let $F$ be a differentially closed field, and let $V\subseteq F^{4n}$ be a strongly $J$-broad variety defined over the field of constants $C$. Then there is a differential field extension $K\supseteq F$ such that $V(K)$ contains an $E_j$-point none of the coordinates of which is constant in $K$. In particular, if $F$ is differentially closed then $V(F)$ contains an $E_J$-point with no constant coordinates.
\end{theorem}

As mentioned above, the original Ax-Schanuel theorem is about the exponential differential equation and the appropriate EC statement in that setting was established by Kirby in \cite{Kirby-semiab}. However, he uses the proof of Ax's theorem which is differential algebraic and allows one to work with an abstract differential field and differential forms. On the other hand, there is no known differential algebraic proof for the Ax-Schanuel theorem for the $j$-function, hence we cannot use its proof (the proof of \cite{Pila-Tsim-Ax-j} uses o-minimality and point counting). Instead we use the statement of the Ax-Schanuel theorem itself, regarding it as a black box without looking inside it. This method is quite general and would probably adapt to other settings to prove an appropriate EC statement as long as there is an appropriate Ax-Schanuel statement. To demonstrate this, in Section \ref{A} we give a proof of EC for the exponential differential equation. In fact, we prove EC for the exponential differential equation in fields with several commuting derivations, which is technically a new result as it is not addressed in Kirby's work. An analogous result for the $j$-function holds as well, with a similar proof. 

Note that Kirby's work on EC for the exponential differential equation, as well as the current paper on EC for the differential equation of the $j$-function, have been motivated by Zilber's work on pseudo-exponential fields and, in particular, his \emph{Exponential Closedness} conjecture for the complex exponential function (see \cite{Zilb-pseudoexp,Bays-Kirby-exp}). Indeed, the differential EC for $\exp$ is a functional/differential analogue of Zilber's conjecture, and the differential EC for $j$ is again its functional/differential analogue for the $j$-function. Naturally, there is a complex existential closedness conjecture for the $j$-function which states roughly that $j$-broad varieties contain points of the form $(\bar{z},j(\bar{z}))$ where each $z_i$ is a complex number in the upper half-plane, and $j$ is the $j$-function. There is also a similar conjecture for $j$ and its derivatives. Both of those are analogues of Zilber's conjecture mentioned above. Eterovi\'c and Herrero have recently made progress towards the complex EC for the $j$-function \cite{eterovic-herrero}.

We remark that the EC theorem for the $j$-function has several important consequences. In particular, in \cite[\S 4]{Aslanyan-adequate-predim} Aslanyan studied $E_j$-reducts of differentially closed fields, that is, reducts of the form $(K;+,\cdot, E_j)$ where $E_j$ is a binary predicate for the solutions of the differential equation of the $j$-function and $K$ is a differentially closed field with a single derivation. He gave a candidate axiomatisation of the common complete theory of those reducts, and proving EC shows that it is indeed an axiomatisation of that theory. It can be thought of as the first-order theory of the differential equation of the $j$-function. EC for $j$ and its derivatives gives a similar axiomatisation of the $E_J$-reducts of differentially closed fields (see \cite[\S 5]{Aslanyan-adequate-predim}). In the terminology of \cite{Aslanyan-adequate-predim} these results mean that the Ax-Schanuel inequality for $E_J$ (as well as $E_j$) is \emph{adequate}.

Using the methods and results of \cite{Aslanyan-adequate-predim} and assuming EC for $E_j$, Aslanyan also  proved a dichotomy result for strongly minimal sets in $E_j$-reducts of differentially closed fields. That statement now becomes unconditional due to this paper (see \cite{Aslanyan-SM-reducts}).

Further, EC (with derivatives) was used in \cite{Aslanyan-weakMZPD} to establish a Zilber-Pink type statement for the $j$-function and its derivatives. Note that Zilber-Pink is a diophantine conjecture generalising Mordell-Lang and Andr\'e-Oort (see, for example, \cite{Zannier-book-unlikely,Habegger-Pila-o-min-certain,Daw-Ren} for details). Two such statements were considered in \cite{Aslanyan-weakMZPD} one of which was proven unconditionally and the other was proven conditionally upon the EC conjecture, so proving the latter here makes it unconditional. 

The current paper may be seen as a continuation of \cite{Aslanyan-adequate-predim}. Nevertheless, the results of this paper, though motivated by \cite{Aslanyan-adequate-predim}, are quite independent from the latter. We have tried to make this paper as self-contained as possible. Most proofs are based on basic facts from differential algebra and algebraic geometry, and of course the Ax-Schanuel theorem.

\section{Ax-Schanuel for the $j$-function}\label{section-ax-sch}

\subsection{Modular polynomials}

The $j$-function is a modular function defined and holomorphic on the upper half-plane $\h:= \{ z \in \C: \im z > 0 \}$. It is invariant under the action of $\SL_2(\mathbb{Z})$ on $\h$, and satisfies certain algebraic ``functional equations'' under the action of $\GL_2^+(\Q)$ -- the group of $2 \times 2$ rational matrices with positive determinant. More precisely, there is a countable collection of irreducible polynomials $\Phi_N\in \Z[X,Y]~ (N\geq 1)$, called \emph{modular polynomials}, such that for any $z_1, z_2\in \h$
$$\Phi_N(j(z_1),j(z_2))=0 \mbox{ for some } N \mbox{ iff } z_2=gz_1 \mbox{ for some } g\in \GL_2^+(\Q). $$

%Moreover, by means of $j$ the quotient $Y':=\SL_2(\mathbb{Z}) \setminus \h$ is identified with $\C$ (thus, $j$ is a bijection from any fundamental domain of $\SL_2(\mathbb{Z})$ to $\mathbb{C}$). 

%In particular, if $\tau\in \h$ is a quadratic number then $j(\tau)$ is algebraic. These numbers are known as \emph{special values} of $j$ or as \emph{singular moduli}. 
%For $w=j(z)$ the image of the $\GL_2^+(\mathbb{Q})$-orbit of $z$ under $j$ is called the \emph{Hecke orbit} of $w$. It obviously consists of the union of solutions of the equations $\Phi_N(X,w)=0,~ N\geq 1$. 
Two elements $w_1,w_2 \in \mathbb{C}$ are called \emph{modularly independent} if they do not satisfy any modular relation $\Phi_N(w_1,w_2)=0$. We refer the reader to \cite{Lang-elliptic} for further details.

\subsection{Differential equation}

The $j$-function satisfies an order $3$ algebraic differential equation over $\mathbb{Q}$, and none of lower order (see \cite{mahler}). Namely, $\Psi(j,j',j'',j''')=0$ where 
$$\Psi(y_0,y_1,y_2,y_3)=\frac{y_3}{y_1}-\frac{3}{2}\left( \frac{y_2}{y_1} \right)^2 + \frac{y_0^2-1968y_0+2654208}{2y_0^2(y_0-1728)^2}\cdot y_1^2.$$

Thus
$$\Psi(y,y',y'',y''')=Sy+R(y)(y')^2,$$
where $S$ denotes the \emph{Schwarzian derivative} defined by $Sy = \frac{y'''}{y'} - \frac{3}{2} \left( \frac{y''}{y'} \right) ^2$ and $R(y)=\frac{y^2-1968y+2654208}{2y^2(y-1728)^2}$ is a rational function.

It is well known that all functions $j(gz)$ with $g \in \SL_2(\mathbb{C})$ satisfy the differential equation $\Psi(y,y',y'',y''')=0$ and all solutions (not necessarily defined on $\mathbb{H}$) are of that form (see \cite{Freitag-Scanlon,Aslanyan-adequate-predim}). %If we allow functions not necessarily defined on $\mathbb{H}$, then all solutions will be of the form $j(gz)$ where $g \in \SL_2(\mathbb{C})$.

Note that for non-constant $y$, the relation $\Psi(y,y',y'',y''') = 0$ is equivalent to $y''' = \eta(y,y',y'')$ where $$\eta(y,y',y'') := \frac{3}{2}\cdot \frac{(y'')^2}{y'} - R(y) \cdot (y')^3$$ is a rational function over $\mathbb{Q}.$

\subsection{Ax-Schanuel}\label{sect-ax-sch}

From now on, $y', y'', y'''$ will denote some variables/coordinates and not the derivatives of $y$. Derivations of abstract differential fields will not be denoted by $'$.

\begin{definition}
Let $(F; +, \cdot, D_1,\ldots,D_m)$ be a differential field with constant field $C = \bigcap_{k=1}^m \ker D_k$. We define a $4$-ary relation $E_J(x, y, y', y'')$  by
$$ \exists y''' \left[ \Psi(y,y',y'',y''')=0 \wedge \bigwedge_{k=1}^m D_ky=y'D_kx \wedge D_ky'=y''D_kx \wedge D_ky''=y'''D_kx\right].$$
The relation $E_J^{\times}(x,y,y',y'')$ is defined by the formula
$$E_J(x,y,y',y'')\wedge x \notin C \wedge y \notin C \wedge y' \notin C \wedge y'' \notin C.$$
\end{definition}

Pila and Tsimerman have established the following transcendence result for the $j$-function which is an analogue of Ax's theorem for the exponential function (see \cite{Ax}, and also Section \ref{A} for a statement). It will play a crucial role in this paper.

\begin{fact}[Ax-Schanuel for $j$, \cite{Pila-Tsim-Ax-j}]\label{j-chapter-Ax-for-j}
Let $(F;+,\cdot,D_1,\ldots,D_m)$ be a differential field with commuting derivations and with field of constants $C$. Let also $(z_i, j_i, j_i', j_i'') \in E_J^{\times}(F),~ i=1,\ldots,n.$ If the $j_i$'s are pairwise modularly independent then $\td_CC\left(\bar{z},\bar{j},\bar{j}',\bar{j}''\right) \geq 3n+\rk (D_kz_i)_{i,k}.$
\end{fact}

As pointed out in the introduction, this statement can be interpreted as follows. Given a variety $V\subseteq F^{4n}$, if $\dim V$ or the dimension of some projections of $V$ is too small, then $V(F)$ cannot contain an $E_J^{\times}(F)$-point whose $j$-coordinates (i.e. the second $n$ coordinates) are pairwise modularly independent. The definition of $J$-broadness given in the next section is based on this observation. Existential Closedness, which is the main result of this paper, is a dual to Ax-Schanuel. It states that if $V$ is $J$-broad then it contains an $E_J$-point.

Note that we will actually work with varieties defined over arbitrary (not necessarily constant) parameters, and we will look for $E_J$-points in those varieties, rather than $E_J^{\times}$-points. We also prove a theorem about certain varieties containing $E_J^{\times}$-points (Theorem \ref{thm-str-$j$-broad-non-const}). It may seem more appropriate to think of that statement as the dual of Ax-Schanuel, however it is weaker than the full EC and the latter is in fact a dual to a ``relative predimension inequality'' governed by Ax-Schanuel. We refer the reader to \cite{Aslanyan-adequate-predim} for details.

\section{Existential Closedness}\label{section-EC}

\subsection{$J$-broad and $J$-free varieties}\label{sect-broad-free-def}

\begin{definition}
Let $n$ be a positive integer, $k \leq n$ and $\bar{i}:=(i_1,\ldots,i_k)$ with $1\leq i_1 < \ldots < i_k \leq n$. Define the projection map $\pr_{\bar{i}}:K^{n} \rightarrow K^{k}$ by
$$\pr_{\bar{i}}:(x_1,\ldots,x_n)\mapsto (x_{i_1},\ldots,x_{i_k}).$$ Also define $\ppr_{\bar{i}}:K^{4n}\rightarrow K^{4k}$ by
$$\ppr_{\bar{i}}:(\bar{x},\bar{y},\bar{y}',\bar{y}'')\mapsto (\pr_{\bar{i}}\bar{x},\pr_{\bar{i}}\bar{y},\pr_{\bar{i}}\bar{y}',\pr_{\bar{i}}\bar{y}'').$$
%Below $\pr_{\bar{i}}$ should always be understood in the second sense.
\end{definition}

\begin{definition}
Let $K$ be an algebraically closed field. An irreducible algebraic variety $V \subseteq K^{4n}$ is \emph{$J$-broad} if for any $1\leq i_1 < \ldots < i_k \leq n$ we have $\dim \ppr_{\bar{i}} (V) \geq 3k$. We say $V$ is \emph{strongly $J$-broad} if the strict inequality $\dim \ppr_{\bar{i}} (V) > 3k$ holds.
\end{definition}

\begin{definition}
An algebraic variety $V \subseteq K^{4n}$ (with coordinates $\left(\bar{x},\bar{y}, \bar{y}', \bar{y}''\right)$) is $J$-\emph{free} if no coordinate is constant on $V$, and it is not contained in any variety defined by an equation $\Phi_N(y_i,y_k)=0$ for some modular polynomial $\Phi_N$ and some indices $i, k$.
\end{definition}

\subsection{EC for varieties defined over arbitrary parameters}

\begin{lemma}\label{lemma-constant}
Let $(F; +, \cdot, D)$ be a differential field and $(z,j,j',j'')\in E_J(F)$ with $j\neq 0, 1728,~ j',j''\neq 0$. If one of the coordinates of $(z,j,j',j'')$ is constant then so are the others.
\end{lemma}
\begin{proof}
If $Dz =0$ then from the definition of $E_J$ we get $Dj = Dj'=Dj''=0$. If $Dj = 0$ or $Dj' = 0$ then since $j', j'' \neq 0$, we must have $Dz=0$. 

Now assume $Dj'' = 0$ and $Dz \neq 0$. Then $\eta(j,j',j'')=0$. For simplicity we will denote $j''$ by $c$. Then $j' = cz + d$ for some constant $d$, and $j = \frac{c}{2}\cdot z^2 + d\cdot z + e$ for some constant $e$. Substituting this into $\eta(j,j',j'')=0$ we get a non-trivial equation for $z$ over $c,d,e$, which implies that $z$ must be constant.
\end{proof}

\begin{theorem}\label{j-free-ec-2}
Let $(F; +, \cdot, D)$ be a differential field, and let $V\subseteq F^{4n}$ be a $J$-broad and $J$-free variety. Then there is a differential field extension $K\supseteq F$ such that $V(K) \cap E_J(K) \neq \emptyset$.
\end{theorem}

\begin{proof}
This is inspired by the proof of \cite[Theorem 3.10]{Kirby-semiab}. An  important difference is that instead of using Proposition 3.7 of that paper, which is an intermediate statement in the proof of Ax-Schanuel, we use Ax-Schanuel itself.

%When $D=0$, that is, $F$ coincides with its constant field, then any point in $V(F^{\alg})$ will trivially be in $E_J(F^{\alg})$ (we extend $D$ to a derivation of $F^{\alg}$). So we assume $D\neq 0$.
Extending $F$ if necessary we may assume that $D\neq 0$. Let $\bar{v}:=\left(\bar{z}, \bar{j}, \bar{j}', \bar{j}''\right)$ be a generic point of $V$ over $F$, and let $K:= F(\bar{v})^{\alg}$ (we can think of $\bar{v}$ as a point defined in a transcendental extension of $F$). We will show that we can extend $D$ to a derivation on $K$ such that $\bar{v}\in E_{J}(K)$. If $\dim V = 3n+l$ then $\td(K/F) = 3n+l$ (where $l\geq 0$).

Let $\Der(K/C)$ be the $K$-vector space of all derivations $\delta:K\rightarrow K$ which vanish on $C$. Consider the subspace  $$\Der(K/D) := \{ \delta \in \Der(K/C): \delta|_F = \lambda D \mbox{ for some } \lambda \in K \}.$$
For every $\delta \in \Der(K/D)$ there is a unique $\lambda\in K$ such that $\delta|_F = \lambda D$, and we denote it by $\lambda_{\delta}$. The map $\varphi: \delta \mapsto \lambda_{\delta}$ is a linear map $\varphi: \Der(K/D)\rightarrow K$. It is surjective, since for every $\lambda \in K$ the map $\lambda D$ can be extended to a derivation of $K$ over $C$. Moreover, $\ker(\varphi) = \Der(K/F)$, hence $$\dim \Der(K/D) = \dim \Der(K/F) +1 = \td(K/F)+1= 3n+l+1. $$

Consider the sequence of inclusions
$$\Der(K/F) \hookrightarrow \Der(K/D) \hookrightarrow \Der(K/C).$$
The space $\Der(K/C)$ can be identified with the dual space of $\Omega(K/C)$ -- the vector space of differential $1$-forms on $K$ over $C$ -- which gives a sequence of surjections\footnote{More precisely, $\Omega(K/C)$ is the vector space generated by the set of symbols $\{ dx: x \in K \}$ modulo the relations
$d(x+y)=dx+dy,~ d(xy)=xdy+ydx,~ dc=0,~ c\in C.$ The map $d:K \rightarrow \Omega(K/C)$ is called the universal derivation on $K$ over $C$ (see also \cite[\S 2]{Aslanyan-weakMZPD}). Then $\Omega(K/D)$ can be defined as the dual of $\Der(K/D)$. It can also be defined as a quotient of $\Omega(K/C)$ so that $\Der(K/D)$  is identified with the dual of $\Omega(K/D)$.}
$$\Omega(K/C) \twoheadrightarrow \Omega(K/D) \twoheadrightarrow \Omega(K/F).$$ This allows us to consider elements of $\Omega(K/C)$ as elements of $\Omega(K/F)$ by identifying them with their images. The following equalities are straightforward:
\begin{gather*}
    \dim \Der(K/F) = \dim \Omega(K/F) =  \td(K/F) = 3n+l,\\
    \dim \Omega(K/D) = \dim \Der(K/D) = 3n+l+1. 
\end{gather*}

Consider the following differential forms in $\Omega(K/C)$:
\begin{equation*}
    \omega_i := dj_i - j_i'dz_i,~ \omega_i' := dj_i' - j_i''dz_i,~ \omega_i'' := dj_i'' - \eta(j_i,j_i',j_i'')dz_i.
\end{equation*}
Note that since we assumed no coordinate is constant on $V$, none of the coordinates of $\bar{v}$ can be in $F$. Therefore $\eta(j_i,j_i',j_i'')$ is well defined.\footnote{ {This is the only place in the proof where we use the fact that no coordinate is constant on $V$. So in fact assuming $y_i\neq 0, 1728$ and $y_i'\neq 0$ on $V$ would suffice.}}

Set $\Lambda(K/C):= \Span_K \{ \omega_i, \omega_i', \omega_i'': i=1,\ldots, n \}\subseteq \Omega(K/C)$ and let $\Lambda(K/D) \subseteq \Omega(K/D)$ and $\Lambda(K/F) \subseteq \Omega(K/F)$ be the images of $\Lambda$ under the canonical surjections.

\begin{claim}
The forms $\omega_i, \omega_i', \omega_i'',~ i=1,\ldots, n$ are $K$-linearly independent in $\Omega(K/F)$, that is, $\dim \Lambda(K/F)  = 3n$.
\end{claim}
\begin{proof}
We proceed by contradiction, so assume $\dim \Lambda(K/F)  < 3n$. Consider the annihilator $\Ann(\Lambda(K/F)) \subseteq \Der(K/F).$ Clearly, $$r:=\dim \Ann(\Lambda(K/F)) = \dim \Omega(K/F) - \dim \Lambda(K/F) > l.$$
It is easy to see that $\Ann(\Lambda(K/F))$ is closed under the Lie bracket, hence we can choose a commuting basis of derivations $D_1,\ldots,D_r \in \Ann(\Lambda(K/F))$ (see \cite[Chapter 0, $\S 5$, Proposition 6]{Kolchin-diff-alg-gp} or \cite[Lemma 2.2]{Singer-noncommuting}). Let $L:= \bigcap_{i=1}^r\ker D_i$ be the field of constants. Then $F \subseteq L \subsetneq K$. 

Since $r>l \geq 0$, at least one of the coordinates of $\bar{v}$ is not in $L$. Let $v_i:=(z_i,j_i,j_i',j_i'')$. By Lemma \ref{lemma-constant} we may assume that for some $t\geq 1$ none of the coordinates of $v_1, \ldots, v_t$ is in $L$ and all coordinates of $v_{t+1},\ldots, v_n$ are in $L$. Let $\bar{u}:= (v_1, \ldots, v_t),~ \bar{w} := (v_{t+1},\ldots,v_n).$

Since $\bar{v}$ is generic in $V$ over $F$, and $V$ is $J$-free, $j_1,\ldots,j_t$ are pairwise modularly independent. By the Ax-Schanuel theorem for several commuting derivations $$\td(L(\bar{u})/L) \geq 3t+\rk(D_i z_k)_{1\leq k\leq t,1\leq i\leq r} = 3t+r.$$ 
 {Here we used the fact that $\rk(D_i z_k)_{1\leq k\leq t,1\leq i\leq r} = r$. To show this we observe that the Jacobian of $\bar{v}$ with respect to $D_1,\ldots,D_r$ has rank $r$ (see, for example, \cite[$\S 2.1$, Claim 1]{Aslanyan-weakMZPD}) and for $k>t$ all coordinates of $v_k$ are in $L$ and hence do not contribute to the rank of the Jacobian. Furthermore, the rows of the Jacobian of $\bar{v}$ corresponding to $j_i,j_i',j_i''$ are linearly dependent on the row corresponding to $z_i$, so those rows do not contribute to the rank either. Thus, $\rk \Jac(\bar{v}) = \rk(D_i z_k)_{1\leq k\leq t,1\leq i\leq r}$.}

Further, $$\td(L/F)\geq \td(F(\bar{w})/F)\geq 3(n-t)$$ for $V$ is $J$-broad. Combining these two inequalities we get $$\td(K/F) = \td(K/L) + \td(L/F) \geq 3t+r + 3(n-t) = 3n+r > 3n+l,$$ which is a contradiction.
\end{proof}

Now we have $\dim \Ann \Lambda(K/D) = \dim \Omega(K/D) - \dim \Lambda(K/D) = 3n+l+1-3n =l+1$ and $\dim \Ann \Lambda(K/F) = \dim \Omega(K/F) - \dim \Lambda(K/F) = l$. Choose a derivation $\delta \in \Ann \Lambda(K/D) \setminus \Ann \Lambda(K/F)$. Then $\delta|_F = \lambda_{\delta} \cdot D$ for some $\lambda \in K$. On the other hand,  $\delta \notin \Ann(\Lambda(K/F))$, therefore $\delta|_F \neq 0$ and $\lambda_{\delta} \neq 0$. Replacing $\delta$ by $\lambda_{\delta}^{-1}\cdot \delta$ we may assume that $\lambda_{\delta} = 1$ and $\delta$ is an extension of $D$ to $K$. 
\end{proof}

Now we are ready to prove Theorem \ref{thm-intro-ec-J}.

\begin{proof}[Proof of Theorem \ref{thm-intro-ec-J}]
If the field $F$ is differentially closed, then it is existentially closed in $K$, that is, a system of differential equations with parameters from $F$ has a solution in $K$ if and only if it has a solution in $F$. Therefore, in this case $V(F)\cap E_J(F)\neq \emptyset$.

{Now we show that in Theorem \ref{j-free-ec-2} the assumption of $J$-freeness may be dropped (cf. \cite[Lemmas 4.43 and 5.24]{Aslanyan-adequate-predim}). We proceed by induction on $n$, the case $n=1$ being trivial. Let $V\seq F^{4n}$ be a $J$-broad variety. If $V$ is $J$-free, then we are done by Theorem \ref{j-free-ec-2}. So we assume $V$ is not $J$-free.} 

If $V$ has a constant coordinate, then assume without loss of generality that one of the coordinates $x_1,y_1,y_1',y_1''$, say $x_1$, is constant, and denote it by $a$. By $J$-broadness, $\dim \ppr_{1}V=3$. Choose elements $b,b',b''$ in an extension $F_1$ of $F$ with $(a,b,b',b'')\in E_J(F_1)$ and let $W \seq F_1^{4(n-1)}$ be the fibre of $V$ above $(a,b,b',b'')$. Then, by the fibre dimension theorem, $W$ is $J$-broad so by the induction hypothesis it contains an $E_J$-point.

{Now assume a modular relation $\Phi_N(y_1,y_2)=0$ holds on $V$. Choose generic constants $a,b,c,d$ in an extension $F_1\supseteq F$ and consider the variety $S\seq F_1^{4n}$ given by the equations $\Phi_N(y_1,y_2)=0,~ \frac{ax_1+b}{cx_1+d}=x_2$ and two more equations obtained by differentiating $\Phi_N(y_1,y_2)=0$ (see \cite[$\S 6.1$]{Aslanyan-weakMZPD}). Let $\bar{i}:=(2,\ldots,n)$ and let $W:=\ppr_{\bar{i}}(V\cap S)\seq F_1^{4(n-1)}$.  {We claim that $W$ is $J$-broad. To this end let $\bar{l}:=(l_1,\ldots,l_m)$ with $2\leq l_1 < \ldots < l_m \leq n$, and let $\hat{l}:=(1,l_1,\ldots,l_m)$. Since $V$ is $J$-broad, $\dim \ppr_{\hat{l}}V \geq 3(m+1)$. Since on $S$ the coordinates $x_2,y_2,y_2',y_2''$ are algebraically related to the coordinates $x_1,y_1,y_1',y_1''$ respectively and there are no relations between any other coordinates, $\dim \ppr_{\hat{l}} S = \dim \ppr_{\bar{l}} S = 4m$, and similarly $$ \dim \ppr_{\bar{l}} (V \cap S) = \dim \ppr_{\hat{l}} (V \cap S) = \dim (\ppr_{\hat{l}} V \cap \ppr_{\hat{l}} S).$$ Using the theorem on dimension of intersection of two varieties (see \cite[\S 1.6.2, Theorem 1.24]{Shafarevich}) and taking into account the fact that the relation $\Phi_N(y_1,y_2)=0$ holds both on $V$ and $S$, we conclude that $$\dim \ppr_{\bar{l}} W =  \dim (\ppr_{\hat{l}} V \cap \ppr_{\hat{l}} S) \geq  \dim \ppr_{\hat{l}} V + 4m+1 - 4(m+1) \geq 3m.$$} 
So by the induction hypothesis $W$ contains an $E_J$-point in a differential field extension of $F_1$. Then the equations defining $S$ allow us to extend it to an $E_J$-point of $V$.} 
%Theorem \ref{thm-intro-ec-J} now follows from Theorem \ref{j-free-ec-2} since it is proven in \cite{Aslanyan-adequate-predim} that in EC one may assume the variety $V$ is $J$-free (see Lemmas 4.43 and 5.24 of that paper).By ,
\end{proof}

\subsection{EC for varieties defined over the constants}

\begin{theorem}\label{thm-str-$j$-broad-non-const}
Let $(F; +, \cdot, D)$ be a differential field, and let $V\subseteq F^{4n}$ be a strongly $J$-broad and $J$-free variety defined over the field of constants $C$. Then there exists a differential field extension $K$ of $F$ such that $V(K) \cap E_J^{\times}(K) \neq \emptyset$. In particular, when $F$ is differentially closed, we have $V(F) \cap E_J^{\times}(F) \neq \emptyset$.
\end{theorem}

This follows from EC by the method of intersecting $V$ with generic hyperplanes (see \cite[Lemma 4.31]{Aslanyan-adequate-predim}). Still, we give a direct proof below, especially as it plays a key role in the proof of a Zilber-Pink type theorem in \cite{Aslanyan-weakMZPD}.

\begin{proof}
We will make use of vector spaces of derivations and differential forms defined in the proof of Theorem \ref{j-free-ec-2}.

Let $\bar{v}:=\left(\bar{z}, \bar{j}, \bar{j}', \bar{j}''\right)$ be a generic point of $V$ over $C$, and let $L:= C(\bar{v})^{\alg}$. If $\dim V = 3n+l$ then $\td(L/C) = 3n+l$ (where $l\geq 1$). So
$$r:=\dim \Ann \Lambda(L/C) = \dim \Omega(L/C) - \dim \Lambda(L/C) \geq l.$$
Note that we actually know that $\dim \Lambda(L/C) = 3n$ and the above inequality is actually an equality, but we do not need it. Pick a commuting basis of derivations $D_1,\ldots,D_r$ of $\Ann \Lambda(L/C)$. Consider the differential field $(L;+,\cdot,D_1,\ldots,D_r)$ and let $C_L := \bigcap_{i=1}^r \ker D_i$ be its constant field. 

 {We claim that none of the coordinates of $\bar{v}$ is in $C_L$. To this end assume by Lemma \ref{lemma-constant} that for some $1\leq t < n$ none of the coordinates of $v_1, \ldots, v_t$ is in $C_L$ and all coordinates of $v_{t+1},\ldots, v_n$ are in $C_L$, where $v_i:=(z_i,j_i,j_i',j_i'')$. Let $\bar{u}:= (v_1, \ldots, v_t),~ \bar{w} := (v_{t+1},\ldots,v_n).$ Since $V$ is strongly $J$-broad, $\td(C_L/C)\geq \td(C(\bar{w})/C)\geq 3(n-t)+1$. Further, by Ax-Schanuel $\td(C_L(\bar{u})/C_L)\geq 3t+r$. Combining these inequalities we get $$\td(L/C) = \td(L/C_L)+\td(C_L/C)\geq 3n+r+1\geq 3n+l+1$$ which is a contradiction. Thus, $t=n$ and $\bar{v}\in (L\setminus C_L)^{4n}.$} 

It is easy to see that we can find a linear combination $\delta$ of $D_1,\ldots,D_r$ such that none of the coordinates of $\bar{v}$ is in $\ker \delta$. Note that $\delta \in \Ann \Lambda(L/C)$. Now in the differential field $(L;+, \cdot, \delta)$ we have $\bar{v}\in E_J^{\times}(L)$. Let $K$ be a common differential field extension of $L$ and $F$ over $C$, that is, $L$ and $F$ can be embedded into $K$ so that the image of $C$ under those embeddings is the same, and the appropriate diagram commutes. Such a field $K$ exists due to the amalgamation property of differential fields. Since $\bar{v} \in V(L) \cap E_J^{\times}(L)$, we also have $\bar{v} \in V(K) \cap E_J^{\times}(K)$.
\end{proof}

\subsection{Generic points}

Often it is important to know that certain varieties contain generic $E_J$-points. In this section we state two such results.

\begin{theorem}\label{thm-gen}
Let $(F; +, \cdot, D)$ be an $\aleph_0$-saturated differentially closed field, and let $V\subseteq F^{4n}$ be a $J$-broad variety defined over a finitely generated subfield $A \subseteq F$. Then $V(F) \cap E_J(F)$ contains a point generic in $V$ over $A$.
\end{theorem}

\begin{theorem}\label{thm-gen-const}
Let $(F; +, \cdot, D)$ be an $\aleph_0$-saturated differentially closed field with field of constants $C$, and let $V\subseteq F^{4n}$ be a strongly $J$-broad variety defined over a finitely generated subfield $C_0 \subseteq C$. Then $V(F) \cap E_J^{\times}(F)$ contains a point generic in $V$ over $C_0$. Furthermore, if $V$ is also $J$-free then $V(F) \cap E_J^{\times}(F)$ contains a point generic in $V$ over $C$.
\end{theorem}

Theorem \ref{thm-gen} follows from Theorem \ref{thm-intro-ec-J} and \cite[Proposition 4.35]{Aslanyan-adequate-predim},  {which is based on Rabinovich's trick to show that a Zariski open subset of an irreducible $J$-broad variety is isomorphic to a $J$-broad variety in a higher dimensional space and the latter contains an $E_J$-point}. The same argument can be applied to deduce the first part of Theorem \ref{thm-gen-const} from Theorem \ref{thm-intro-ec-J-nonconst}. The second part of Theorem \ref{thm-gen-const} follows from \cite[Proposition 4.29]{Aslanyan-adequate-predim}.

\subsection{EC for $j$ without derivatives}\label{sect-EC-j-noder}

In a differential field we define a binary relation $E_j(x,y)$ as the projection of $E_J$ onto the first two coordinates, in other words we ignore the derivatives. The theorems proved above obviously imply similar statements for $E_j$, in particular, Theorem \ref{thm-intro-ec-j}. First we define $j$-broad varieties.

\begin{definition}
\begin{itemize}[leftmargin=0.7cm]
    \item[] 
    
    \item For positive integers $n\geq k$ and $1\leq i_1 < \ldots < i_k \leq n$ define a projection map $\pi_{\bar{i}}:K^{2n} \rar K^{2k}$ by
    $\pi_{\bar{i}}:(\bar{x},\bar{y})\mapsto (\pr_{\bar{i}}\bar{x},\pr_{\bar{i}}\bar{y}).$

    \item Let $K$ be an algebraically closed field. An irreducible algebraic variety $V \subseteq K^{2n}$ is \emph{$j$-broad} if for any $1\leq i_1 < \ldots < i_k \leq n$ we have $\dim \pi_{\bar{i}} (V) \geq k$.
\end{itemize}
\end{definition}

\begin{proof}[Proof of Theorem \ref{thm-intro-ec-j}]
If $V \subseteq K^{2n}$ is $j$-broad then $\tilde{V}:=V \times K^{2n}$ is obviously $J$-broad. Hence, by Theorem \ref{thm-intro-ec-J} there is a differential field extension $K\supseteq F$ and a point $(\bar{z},\bar{j},\bar{j}',\bar{j}'')\in \tilde{V}(K)\cap E_J(K)$. Then $(\bar{z},\bar{j})\in V(K)\cap E_j(K)$.
\end{proof}

\section{EC in fields with several commuting derivations}\label{A}

In this section we give a proof of EC for the exponential differential equation in fields with several commuting derivations. Strictly speaking, it is a new result since Kirby addressed the EC question only for ordinary differential fields. Note that even in the case of a single derivation, our proof differs from Kirby's proof (see \cite[Theorem 3.10]{Kirby-semiab}) in that we use Ax-Schanuel, rather than its proof. The proof, which is a slight generalisation of the arguments presented above, also works for the $j$-function and establishes an EC result for $E_J$ in fields with several commuting derivations. Moreover, we believe that our idea of using Ax-Schanuel to prove EC is quite general and will go through in other settings too, and proving EC for the exponential differential equation in this section supports this speculation.

The EC statement for the exponential differential equation presented below differs from the EC statement of \cite[Definition 2.30]{Kirby-semiab}. In fact, the latter is not completely correct, and it was corrected in \cite[Definition 10.3]{Bays-Kirby-exp} (see also \cite[Fact 2.3]{Kirby-blurred}), which is the statement we consider in this section.

In a differential field $(F; +, \cdot, D_1,\ldots,D_m)$ we define a binary relation $\Exp(x,y)$ by the formula $\bigwedge_{k=1}^m D_k y = y D_k x.$

\begin{fact}[Ax-Schanuel, \cite{Ax}]\label{Ax-Schanuel-exp}
Let $(F;+,\cdot,D_1,\ldots,D_m)$ be a differential field with field of constants $C = \bigcap_{k=1}^m \ker D_k$. Let also $(x_i, y_i) \in \Exp(F),~ i=1,\ldots,n,$ be such that $x_1,\ldots,x_n$ are $\Q$-linearly independent mod $C$, that is, they are $\Q$-linearly independent in the quotient vector space $F/C$.
Then
$\td_CC\left(\bar{x},\bar{y}\right) \geq n+\rk (D_kx_i)_{i,k}.$
\end{fact}

For a field $F$ let $\G_{\aaa}(F)$ and $\G_{\m}(F)$ denote the additive and multiplicative groups of $F$ respectively, and for a positive integer $n$ let $G_n:=\G^n_{\aaa} \times \G^n_{\m}$. For a $k \times n$ matrix $M$ of integers we define $[M]:G_n(F) \rightarrow G_k(F)$ to be the map given by $[M]:(\bar{x},\bar{y}) \mapsto (u_1,\ldots,u_k, v_1,\ldots, v_k)$ where
$$u_i = \sum_{j=1}^n m_{ij}x_j \mbox{ and } v_i = \prod_{j=1}^n y_j^{m_{ij}}.$$

\begin{definition}\label{rotund}
An irreducible variety $V \subseteq G_n(F)$ is \emph{rotund} if for any $1 \leq k \leq n$ and any $k\times n$ matrix $M$ of integers $\dim [M](V) \geq \rk M$. %If for any non-zero $M$ the stronger inequality $\dim [M](V) \geq \rk M + 1$ holds then we say $V$ is \emph{strongly rotund}.
\end{definition}

%\begin{remark}
%$J$-broadness of a variety is the analogue of rotundity for the $j$-function.
%\end{remark}

\begin{theorem}
Let $(F; +, \cdot, D_1, \ldots, D_m)$ be a differential field with $m$ commuting derivations, and let $V\subseteq F^{2n}$ be a rotund variety. Then there exists a differential field extension $K$ of $F$ such that $V(K) \cap \Exp(K) \neq \emptyset$. In particular, when $F$ is differentially closed, $V(F) \cap \Exp(F) \neq \emptyset$.
\end{theorem}

\begin{proof}
Let $\bar{v}:=(\bar{x}, \bar{y})$ be a generic point of $V$ over $F$. Set $K:= F(\bar{v})^{\alg}$. If $\dim V = n+l$ for some $l\geq 0$ then $\td(K/F) = n+l$. 

We may assume the derivations $D_1, \ldots, D_m$ are $F$-linearly independent. Let $\Delta:= \Span_F\{ D_1, \ldots, D_m \} \subseteq \Der(F/C)$.
Consider the space of derivations $$\Der(K/\Delta) := \{ D \in \Der(K/C): D|_F = \lambda_1 D_1 + \ldots + \lambda_m D_m \mbox{ with } \lambda_i \in K \}.$$
Extending $F$ if necessary, we may assume that for each $i$ there is $t_i \in F$ such that $D_it_i\neq 0$ and $D_kt_i =0$ whenever $k\neq i$. Therefore, for every $D \in \Der(K/\Delta)$ there is a unique tuple $\bar{\lambda}\in K^m$ such that $D|_F = \lambda_1 D_1 + \ldots + \lambda_m D_m$, and we denote it by $\bar{\lambda}_D$. The map $\varphi: D \mapsto \bar{\lambda}_D$ is a linear map $\varphi:\Der(K/\Delta)\rightarrow K^m$. It is clearly surjective and $\ker(\varphi) = \Der(K/F)$, hence $$\dim \Der(K/\Delta) = \dim \Der(K/F) + m = n+l+m. $$

As in the proof of Theorem \ref{j-free-ec-2}, we have a sequence of inclusions
$$\Der(K/F) \hookrightarrow \Der(K/\Delta) \hookrightarrow \Der(K/C),$$
and a dual sequence of surjections
$$\Omega(K/C) \twoheadrightarrow \Omega(K/\Delta) \twoheadrightarrow \Omega(K/F).$$ 

Consider the differential forms $\omega_i := dy_i - y_idx_i\in \Omega(K/C)$. Denote $\Lambda(K/C):= \Span_K \{ \omega_i: i=1,\ldots, n \}\subseteq \Omega(K/C)$ and let $\Lambda(K/\Delta) \subseteq \Omega(K/\Delta)$ and $\Lambda(K/F) \subseteq \Omega(K/F)$ be the images of $\Lambda(K/C)$ under the canonical surjections.

\begin{claim}
The forms $\omega_i,~ i=1,\ldots, n,$ are $K$-linearly independent in $\Omega(K/F)$, that is, $\dim \Lambda(K/F)  = n$.
\end{claim}
\begin{proof}
Assume $\dim \Lambda(K/F)  < n$. Consider the annihilator $\Ann(\Lambda(K/F))$ as a subspace of $\Der(K/F).$ It is clear that $$r:=\dim \Ann(\Lambda(K/F)) = \dim \Omega(K/F) - \dim \Lambda(K/F) > l.$$
The space $\Ann(\Lambda(K/F))$ is closed under the Lie bracket, hence we can choose a commuting basis of derivations $\delta_1,\ldots,\delta_r \in \Ann(\Lambda(K/F))$. Then $\delta_i y_k = y_k \delta_i x_k$ for all $i, k$. Let $L:= \bigcap_{i=1}^r\ker \delta_i$ be the field of constants. Then $F \subseteq L \subsetneq K$. 

We choose a maximal subtuple of $\bar{x}$ which is $\Q$-linearly independent modulo $L$. Assume without loss of generality that $(x_1,\ldots,x_t)$ is such a subtuple. Observe that $t>0$ since otherwise $L=K$. For every $i>t$ there are integers $\beta_i \neq 0,~ \alpha_1^i, \ldots, \alpha_t^i$ and an element $l_i \in L$ such that $$\alpha_1^ix_1 + \ldots + \alpha_t^ix_t + \beta_i x_i = l_i.$$

Let $u_i:= y_1^{\alpha_1^i}\cdots y_t^{\alpha_t^i} \cdot y_i^{\beta_i}$. Obviously $\delta_ku_i = u_i \delta_kl_i = 0$ for all $i,k$, hence $u_i \in L$. 

By rotundity of $V$ we have $\td(L/F) \geq \td(F(\bar{l},\bar{u})/F)\geq n-t,$ and by the Ax-Schanuel theorem $\td(L(\bar{x},\bar{y})/L) \geq t+\rk(\delta_ix_k)_{k,i} = t+r.$ Thus $$\td(K/F) = \td(K/L) + \td(L/F) \geq t+r + n-t = n+r>n+l,$$ which is a contradiction.
\end{proof}

Now we have $\dim \Ann \Lambda(K/\Delta) = \dim \Omega(K/\Delta) - \dim \Lambda(K/\Delta) = (n+l+m)-n = l+m$ and $\dim \Ann \Lambda(K/F) = \dim \Omega(K/F) - \dim \Lambda(K/F) = l$. Choose derivations $\delta_1, \ldots, \delta_m \in \Ann \Lambda(K/\Delta) $ which are $K$-linearly independent modulo $\Ann \Lambda(K/F)$. Each $\delta_i|_F$ is a $K$-linear combination of $D_k$'s, so in a matrix form $\delta|_F = M \cdot D,$ where $\delta$ and $D$ are the column vectors of the derivations $\delta_i$ and $D_i$ respectively and  $M\in \GL_m(K)$ (it is invertible for otherwise a $K$-linear combination of $\delta_i$'s would be in $\Der(K/F)$). Replacing $\delta$ with $M^{-1} \cdot \delta$ we may assume $M$ is the identity matrix and hence each $\delta_i$ is an extension of $D_i$. 
%For the second part of the theorem we appeal to existential closedness of $F$ in $K$.
\end{proof}

\iffalse
\begin{remark}
Notice that we did not need to assume that the variety $V$ is free here, while in the case of $j$ it seems our argument would fail for non-free varieties. Of course, as we have pointed out, in both settings EC can be reduced to free varieties.
\end{remark}
\fi

The following theorem for the $j$-function can be proven similarly.

\begin{theorem}
Let $F$ be a differential field with several commuting derivations, and let $V\subseteq F^{4n}$ be a $J$-broad variety. Then there exists a differential field extension $K$ of $F$ such that $V(K) \cap E_J(K) \neq \emptyset$. In particular, when $F$ is differentially closed, $V(F) \cap E_J(F) \neq \emptyset$.
\end{theorem}

\subsection*{Acknowledgements.} We are grateful to the referee for useful comments that helped us improve the presentation of the paper.

%\addcontentsline {toc} {section} {Bibliography}
\bibliographystyle {alpha}
\bibliography {ref}

\end{document}